# MAXIMAL NILPOTENT QUOTIENTS OF 3-MANIFOLD GROUPS

PETER TEICHNER

ABSTRACT. We show that if the lower central series of the fundamental group of a closed oriented 3-manifold stabilizes then the maximal nilpotent quotient is a cyclic group, a quaternion 2-group cross an odd order cyclic group, or a Heisenberg group. These groups are well known to be precisely the nilpotent fundamental groups of closed oriented 3-manifolds.

## 1. INTRODUCTION

There are many different approaches to the study of 3-manifolds. There is a flourishing combinatorial school from Dehn, Papakyriaskopolous, Haken and Waldhausen to Gordon and Luecke, which has shown that many 3-manifold questions can be reduced to questions about the fundamental group. A co-equal off-shoot of the combinatorial school folds in dynamics and complex analysis. This is Thurston's program on geometrizing certain *characteristic* and *simple* pieces of 3-manifolds by showing that each is modelled locally on one of the eight 3-dimensional geometries (of which only the hyperbolic case is not fully understood). A third perspective on 3-manifolds is through quantuum field theory. The ideas of Witten, Jones, Vassiliev and many others have inspired tremendous activity and, in time, may contribute substantially to the topological understanding of 3-manifolds.

This paper takes a fourth perspective by looking at a 3-manifold through nilpotent eyes, observing only the tower of nilpotent quotients of the fundamental group, but never the group itself. This point of view has a long history in the study of link complements and it arises naturally if one studies 3- and 4-dimensional manifolds together. For example, Stallings proved that for a link in $S^3$ certain nilpotent quotients of the fundamental group of the link complement are invariants of the topological concordance class of the link. These quotients contain the same information as Milnor's $\bar{\mu}$-invariants which are generalized linking numbers. For precise references about this area of research and the most recent applications to 4-manifolds see [4].

Turaev [9] seems to have been the first to consider nilpotent quotients of closed 3-manifold rather than link complements. Much earlier, the nilpotent fundamental groups of closed 3-manifolds were classified. Thomas [8] showed in particular that statements (1) and (2) in the following Theorem 1 are equivalent.

---

Research at MSRI was supported by a fellowship from the Miller foundation.

**Theorem 1.** *For a nilpotent group $N$ the following statements are equivalent:*

(1) $N$ *is a (finite or infinite) cyclic group* $\mathbb{Z}/n$*, a product* $Q_{2^n} \times \mathbb{Z}/(2k+1)$*, or a Heisenberg group* $H_n$*.*

(2) $N$ *is the fundamental group of a closed orientable 3-manifold.*

(3) $N$ *is finitely generated and there exists a class* $m \in H_3(N)$ *such that the cap-product with* $m$ *induces an epimorphism* $H^1(N) \to H_2(N)$ *and an isomorphism* Torsion $H^2(N) \to$ Torsion $H_1(N)$*.*

(4) $N$ *is the maximal nilpotent quotient of the fundamental group of a closed orientable 3-manifold.*

Recall that a group has a maximal nilpotent quotient if and only if it's lower central series stabilizes. In (1) above the Heisenberg groups $H_n$ are the central extensions of $\mathbb{Z}^2$ by $\mathbb{Z}$ classified by the Euler class $n \in \mathbb{Z} \cong H^2(\mathbb{Z}^2; \mathbb{Z})$. They occur as the fundamental groups of orientable circle bundles over the 2-torus. Euler class $n = 0$ corresponds to the 3-torus. The infinite cyclic group is $\pi_1(S^1 \times S^2)$ and the finite cyclic and generalized quaternion groups $Q_{4k}$ are subgroups of $SU(2)$. (Here $Q_{4k} := \langle x, t \mid txt^{-1} = x^{-1}, x^k = t^2 \rangle$ has order $4k$.) Finally, a product $Q_{2^n} \times \mathbb{Z}/(2k+1)$ can be embedded into $SO(4)$ such that it acts freely on $S^3$, see [10]. Thus these finite groups are fundamental groups of 3-dimensional homogenous spaces.

This shows that (1) implies (2) in the above theorem. The other easy fact is that (2) implies (3): If $N$ is the fundamental group of a closed oriented 3-manifold then $m$ can be taken to be the image of it's fundamental class. From [9, Thm.2] it easily follows that statement (3) implies (4). Thus the aim of this paper is to show that (4) implies (1).

Theorem 1 can be read in two essentially different ways: A group theorist may take (3) as a homological characterization of the family of groups in (1). A 3-manifold topologist would probably prefer the point of view of (4), i.e. that there are very few possibilities for the lower central series of the fundamental group of a closed oriented 3-manifold: Either it descends forever or it stabilizes with the maximal nilpotent quotient being one of the groups in (1).

It is easy to construct examples exhibiting both phenomena: For example, take any surgery description of one of the 3-manifolds described above. Tying a local knot into one of the components does not change the lower central series of the 3-manifold group but it drastically changes the 3-manifold itself. Conversely, if $H_1(M)$ has more than 3 generators then it has no maximal nilpotent quotient: This follows from the fact that $H_1(H_n) \cong \mathbb{Z}^2 \times \mathbb{Z}/n$ and $H_1(Q_{2^n}) \cong \mathbb{Z}/2 \times \mathbb{Z}/2$ for $n > 2$ (and $Q_4 = \mathbb{Z}/4$). Moreover, the only groups in list (1) that have nilpotency class $> 3$ are the one containing $Q_{2^n}, n > 3$. In this single case it might be more difficult to decide whether $\pi_1 M$ has a maximal nilpotent quotient because $Q_{2^n}$ modulo the $k$-th term of the lower central series is the group $Q_{2^k}$ for all $k > 2$.

This problem vanishes if one considers a rational version of Theorem 1. We will prove it in analogy to the integral case. This rational result was first obtained in [3] using completely different methods, namely methods from rational homotopy theory. This approach however does not seem to provide a proof of Theorem 1 but it motivated the research in this paper.

**Theorem 2.** *For a torsionfree nilpotent group $N$ the following statements are equivalent:*

(1) *$N$ is trivial, infinite cyclic or a Heisenberg group.*
(2) *$N$ is the fundamental group of a closed orientable 3-manifold.*
(3) *$N$ is finitely generated and there exists a class $m \in H_3(N; \mathbb{Q})$ such that the cap-product with $m$ induces an epimorphism $H^1(N; \mathbb{Q}) \to H_2(N; \mathbb{Q})$*
(4) *$N$ is the maximal torsionfree nilpotent quotient of the fundamental group of a closed orientable 3-manifold.*

The paper is organized as follows: In Section 2 we give the necessary definitions and proof some (probably well-known) results on nilpotent groups. Section 3 contains the proof of Theorem 1 modulo three Propositions. Our proof is modelled on a proof of Thomas' theorem which we give as a "warm up". Finally, Section 4 contains the proofs of the three Propositions used in Section 3 and the proof of Theorem 2.

## 2. Some nilpotent group theory

The *lower central series* of a group $G$ is defined by $G_1 := G$ and $G_{k+1} := [G, G_k]$ for $k > 1$. $G$ is *nilpotent* if $G_k = 1$ for some $k$ and the smallest such $k$, if it exists, is called the *class* of $G$. Thus abelian groups are precisely the groups of class 2. Any group $G$ has the nilpotent quotients $G/G_k$ and it has a *maximal nilpotent quotient* if and only if $G_k = G_{k+1}$ for some $k$. We define the *rank* of an abelian group $A$ to be the dimension of $A \otimes \mathbb{Q}$. For a nilpotent group $N$ we define it's rank to be the sum of the ranks of the abelian groups $N_k/N_{k+1}$.

**Lemma 1.** *A nilpotent group $N$ is finite if and only if $H_1(N)$ is finite.*

*Proof.* The proof is an induction on the class of $N$. If the class is 2 the statement is true because $N \cong H_1(N)$. Assume that $N/N_k$ is finite. Then $H_2(N/N_k)$ is also a finite group, see [1]. The same reference explains the *5-term exact sequence* for groups: Given a short exact sequence

$$1 \longrightarrow N \longrightarrow G \longrightarrow Q \longrightarrow 1$$

of groups, the bottom part of the corresponding Leray-Serre spectral sequence is an exact sequence of homology groups (with integral coefficients)

$$H_2(G) \longrightarrow H_2(Q) \longrightarrow N/[N, G] \longrightarrow H_1(G) \longrightarrow H_1(Q) \longrightarrow 0.$$

Applying this to the central extension

$$1 \to N_k/N_{k+1} \to N/N_{k+1} \to N/N_k \to 1$$

implies that $N/N_{k+1}$ is also finite.                                    □

**Corollary 2.** *For a nilpotent group $N$, the set of all elements of finite order is a (characteristic) subgroup of $N$.*

*Proof.* Let $g, h$ be elements of finite order in $N$. We need to show that their product $g \cdot h$ is still of finite order. Let $G$ be the subgroup of $N$ generated by $g$ and $h$. It suffices to show that $G$ is finite. Clearly, $H_1(G)$ is finite and since $G$ is nilpotent the result follows from Lemma 1.                                    □

The above subgroup is called the *torsion subgroup* $\mathrm{Tor}(N)$ of $N$. $N$ is *torsionfree* if $\mathrm{Tor}(N) = 1$. In the following lemma we will use a commutator identity which holds in any group $G$. Namely, for elements $a, b \in G$ and $n \in \mathbb{N}$ one has

$$[a^n, b] = [a, b]^{(a^{n-1})} \cdot [a, b]^{(a^{n-2})} \cdot \dots \cdot [a, b]$$

if one uses the conventions $[a, b] := a \cdot b \cdot a^{-1} \cdot b^{-1}$ and $b^a := a \cdot b \cdot a^{-1}$.

**Lemma 3.** *A nilpotent group $N$ is torsionfree if and only if it's center $C$ and $N/C$ are torsionfree.*

*Proof.* The only thing to show is that if $N$ is torsionfree and $x \in N$ satisfies $x^n \in C$ for some $n \in \mathbb{N}$ then $x \in C$. Take any element $g \in N$ and define $x_k$ by $x_1 := [x, g]$ and $x_{k+1} := [x, x_k]$. Since $N$ is nilpotent, we know that $x_{c+1} = 1$ for some $c$. This means that $x$ commutes with $x_c$ and thus the above commutator identity simplifies to give

$$[x^n, x_{c-1}] = [x, x_{c-1}]^n = (x_c)^n$$

But by assumption $x^n$ is central which implies $(x_c)^n = 1$. Since $N$ is torsionfree this indeed shows that $x_c = 1$. Continuing in exactly the same manner leads to $1 = x_{c-1} = \dots = x_1 = [x, g]$. Since $g$ was arbitrarily, we can conclude that $x$ is central.                                    □

The main result of this section follows. We use the notation $\mathrm{hd}_{\mathbb{Z}}(G)$ for the $\mathbb{Z}$-*homological dimension* of a group $G$, i.e. the smallest $n \in \mathbb{N} \cup \infty$ such that $H_i(G; \mathbb{Z})$ vanishes for all $i > n$. We also introduce some notation which will be used in the proof and throughout the rest of the paper:

Let $H$ be a group and $P$ an $H$-module. Then the 0-th homology with twisted coefficients $H_0(H; P) \cong P/I(H) \cdot P$ is isomorphic to the *cofixed point set* $P_H$ of $P$ under the $H$-action, i.e. the largest quotient module of $P$ on which $H$ acts trivially, [1]. Here $I(H)$ is the augmentation ideal of the group ring $\mathbb{Z}[H]$.

Similarly, $H^0(H; P)$ is isomorphic to the *fixed point set* $P^H$ of $P$ under the $H$-action, i.e. the largest submodule of $P$ on which $H$ acts trivially.

**Lemma 4.** *For a finitely generated nilpotent group $N$ the following statements are equivalent:*

(i) $\mathrm{hd}_{\mathbb{Z}}(N) < \infty$.

(ii) $N$ *is torsionfree.*

(iii) $K(N, 1)$ *is homotopy equivalent to a closed orientable manifold of dimension* $\mathrm{rank}(N) = \mathrm{hd}_{\mathbb{Z}}(N)$. *More precisely, $K(N, 1)$ is homotopy equivalent to an iterated circle bundle with structure groups $U(1)$.*

*Proof.* The conclusions (iii) $\Rightarrow$ (i) and (iii) $\Rightarrow$ (ii) are obvious. (ii) $\Rightarrow$ (iii) follows by induction from Lemma 3. One just has to observe that $H^2(N; \mathbb{Z})$ classifies central extensions of $N$ by $\mathbb{Z}$ as well as principal circle bundles over $K(N, 1)$ with structure group $U(1)$. The induction starts with the fact that the $r$-torus is a $K(\mathbb{Z}^r, 1)$.

To show (i) $\Rightarrow$ (ii) we will induct on the rank of $N$. If $\mathrm{rank}(N) = 0$ then $H_1(N)$ is finite and by Lemma 1 $N$ is also finite. But the integral homology groups of a finite group are nontrivial in infinitely many dimensions [7]. So our assumption $\mathrm{hd}_{\mathbb{Z}}(N) < \infty$ implies $N = 1$.

Now assume that $\mathrm{rank}(N) > 0$. Then $N$ and thus $H_1(N)$ are infinite and we get an extension

$$1 \to U \to N \to \mathbb{Z} \to 1$$

of nilpotent groups with $\mathrm{rank}(U) < \mathrm{rank}(N)$. Below we show that $\mathrm{hd}_{\mathbb{Z}}(U) < \infty$ which implies by induction that $U$ and hence $N$ are torsionfree.

Let $t \in N$ be an element which maps to a generator of $\mathbb{Z}$ in the above extension. Then the infinite cyclic group $\langle t \rangle$ acts by conjugation on $U$ and thus on $H_i(U)$. The Wang sequence for the above extension gives exact sequences for any $i > 0$:

$$0 \to H_i(U)_{\langle t \rangle} \to H_i(N) \to H_{i-1}(U)^{\langle t \rangle} \to 0$$

Our claim that $\mathrm{hd}_{\mathbb{Z}}(U) < \infty$ now follows from the following

**Lemma 5.** *The $\langle t \rangle$-modules $H_i(U)$ are nilpotent, i.e. there exists integers $N_i$ such that $(t-1)^{N_i}$ is the zero-map on $H_i(U)$. In particular, if $H_i(U)$ is nontrivial then so are the fixed and cofixed point sets under the $\langle t \rangle$-action.*

To prove the lemma, define $U_k$ by $U_1 := U$ and $U_{k+1} := [N, U_k]$. Since $N$ is nilpotent we have $U_c = 1$ for some $c$ which is smaller or equal to the class of $N$. Moreover, $t$ acts as the identity on the quotients $U_k / U_{k+1}$. By induction, we assume that $\langle t \rangle$ acts nilpotently on the modules $H_i(U/U_k)$. The Leray-Serre spectral sequence for the extension

$$1 \to U_k / U_{k+1} \to U / U_{k+1} \to U / U_k \to 1$$

has $E^2_{p,q}$-terms $H_p(U/U_k; H_q(U_k/U_{k+1}))$ which are then also nilpotent $\langle t \rangle$-modules. Consequently, the $E^{\infty}_{p,q}$-terms are all nilpotent $\langle t \rangle$-modules and so is $H_{p+q}(U/U_{k+1})$ as a (finitely) iterated extension of these modules. $\square$

The following result of B. Dwyer [2] will be essential for our proof.

**Theorem 3.** *Let $N$ be a finitely generated nilpotent group and $P$ a finitely generated $\mathbb{Z}[N]$-module. If $H_0(N; P) = 0$ then $H_i(N; P) = 0$ for all $i \geq 0$.*

One can extend Dwyer's theorem to cohomology groups.

**Corollary 6.** *Let $N$ be a finitely generated nilpotent group and $P$ a finitely generated $\mathbb{Z}[N]$-module. If $H_0(N; P) = 0$ then $H^i(N; P) = 0$ for all $i \geq 0$.*

*Proof.* Let $Z$ be a nontrivial central cyclic subgroup of $N$. It is enough to show that the $E_2$-term of the Serre spectral sequence

$$E_2^{p,q} = H^p(N/Z; H^q(Z; P)) \Longrightarrow H^{p+q}(N; P)$$

vanishes identically. By induction it suffices to show that the assumption of Dwyer's theorem is satisfied for the $N/Z$-modules $H^q(Z; P)$. Since $Z$ is a cyclic group, these twisted cohomology groups are fully understood, see [1, p.58]. In particular, one can show that $H_0(N/Z; H^q(Z; P)) = 0$ by using Lemma 7 below which is an easy consequence of Dwyer's theorem and the long exact coefficient sequence [1, p.71]. $\square$

**Lemma 7.** *Let $N$ be a finitely generated nilpotent group and*

$$0 \to P' \to P \to P'' \to 0$$

*an extension of $\mathbb{Z}[N]$-modules with $P$ (and thus $P', P''$) finitely generated. Then $P_N = 0$ implies that $(P')_N = (P'')_N = 0$.*

## 3. Outline of the proof of Theorem 1

As a warm up, we first give a proof of Thomas' theorem [8] that (2) implies (1) in Theorem 1. So let $N$ be the nilpotent fundamental group of a closed orientable 3-manifold $M$. Then $M$ is the connected sum of prime 3-manifolds $M_i$. Therefore, $N$ is the free product of the fundamental groups $\pi_1 M_i$. Since $N$ is nilpotent, only one of these groups can be nontrivial and thus we may assume that $N$ is the fundamental group of a *prime* 3-manifold $M$. The argument splits now into several cases.

**Case I**: $N$ is a finite group.

Then $N$ acts freely on the homotopy 3-sphere $\widetilde{M}$ and thus has 4-periodic cohomology. Since $N$ is a finite nilpotent group, it is the direct product of it's $p$-Sylow subgroups [5] which are also 4-periodic [1, p.156]. The same reference shows that the $p$-Sylow subgroups are cyclic for $p$ odd, and cyclic or generalized quaternion for $p = 2$. Therefore, $N$ is the direct product of such groups which we wanted to show.

**Case II**: $N$ has infinite order.

This case splits naturally into two subcases depending on whether or not $\pi_2 M$ is trivial. If not, then the prime 3-manifold $M$ is homeomorphic to $S^1 \times S^2$ and hence $N$ is infinite cyclic.

If $\pi_2 M = 0$ then $M$ is a $K(N,1)$ and hence $N$ is torsionfree of rank 3. The center $C$ of $N$ is nontrivial and torsionfree. Moreover, $N/C$ is also torsionfree by Lemma 3. Note that $C = \mathbb{Z}^2$ is impossible: It would imply that $N/C \cong \mathbb{Z}$ is generated by one element which commutes with itself and $C$. Therefore, either $C = N$ or $C \cong \mathbb{Z}$. In both cases $N$ is a central extension of $\mathbb{Z}^2$ by $\mathbb{Z}$ and thus one of the Heisenberg groups.

This finishes the proof that (2) implies (1) in Theorem 1. We now outline the proof that (4) also implies (1) which is the only part of Theorem 1 we need to prove. So let $M$ be a closed orientable 3-manifold whose fundamental group allows a maximal nilpotent quotient $N$. Consider the fibration

$$M' \to M \to K(N,1)$$

induced from the quotient map $\pi_1 M \to N$. Then $M'$ has the homotopy type of the covering space of $M$ for which $N$ is the group of deck transformations. We will do calculations with the Serre spectral sequence

$$E^2_{p,q} = H_p(N; H_q(M')) \Longrightarrow H_{p+q}(M)$$

for the above fibration. Since $N$ is the *maximal* nilpotent quotient of $\pi_1 M$ we get $H_1(M')_N = \pi_1 M'/[\pi_1 M', \pi_1 M] = 0$. Applying Dwyer's Theorem 3 for $P := H_1(M')$ we get $E^2_{i,1} = 0$ for all $i \geq 0$ in our spectral sequence. As above, our argument splits now into several cases.

**Case I**: $N$ is a finite group.

Then $M'$ is up to homotopy a compact 3-manifold and thus $H_3(M') = \mathbb{Z}$. Moreover, we will prove the following proposition in Section 4.

**Proposition 1.** *In the above situation, $H_2(M')_N = 0$.*

This result, Dwyer's Theorem 3 and the fact that the covering map $M' \to M$ has degree $|N|$ imply that we get an element of order $|N|$ in $H_3(N) \cong H^4(N)$. By [1, p.154] $N$ then has 4-periodic cohomology. The proof in Case I now concludes exactly as in the warm up.

**Case II**: $N$ has infinite order.

Then $H_i(M') = 0$ for all $i \geq 3$ and we will prove the following result in Section 4.

**Proposition 2.** *In this situation*

$$H_2(M') \cong \begin{cases} \mathbb{Z} & \text{if } \operatorname{rank}(N) = 1, \\ 0 & \text{else.} \end{cases}$$

Therefore, as in the warm up, we naturally have to consider two cases.

**Case IIa**: $\operatorname{rank}(N) = 1$, i.e. $N/\operatorname{Tor}(N)$ is infinite cyclic.

Surprisingly, it takes quite a bit of work to show that in fact $\operatorname{Tor}(N)$ must be trivial in this case. This will be shown in Proposition 3, Section 4.

**Case IIb**: $\operatorname{rank}(N) > 1$.

Then our calculations above show that the map $M \to K(N, 1)$ is an integral homology isomorphism. Therefore, $N$ is torsionfree of rank 3 by Lemma 4. This implies exactly as in the warm up that $N$ is one of the Heisenberg groups.    □

## 4. Proofs of the Propositions

We have to fill the gaps in the proof of Theorem 1. Recall that $M$ is a closed orientable 3-manifold whose fundamental group allows a maximal nilpotent quotient $N$. We consider the fibration

$$M' \to M \to K(N, 1)$$

induced from the quotient map $\pi_1 M \to N$.

**Case I**: $N$ is a finite group.

Then we need to prove that $H_2(M')_N = 0$ which is the conclusion of Proposition 1.

*Proof.* Set $P := H_1(M')$ which is a finitely generated abelian group. To show $H_2(M')_N = 0$ note that by Poincaré duality $H_2(M') \cong H^1(M') \cong P^*$. From $P_N = 0$ we conclude that $(P^*)^N = 0$ and thus

$$(P^*)_N = H_0(N; P^*) \cong \hat{H}^{-1}(N; P^*) \cong H^1(N; P^{**}),$$

see [1, p.134 and p.148, exercise 3]. Now $P^{**} \cong P/\operatorname{Torsion}$ is a quotient of $P$ and by Lemma 7 it still satisfies $(P^{**})_N = 0$. By Corollary 6 we get $H_2(M')_N \cong (P^*)_N \cong H^1(N; P^{**}) = 0$.    □

**Case II**: $N$ has infinite order.

Now $H_i(M') = 0$ for all $i \geq 3$ and Poincaré duality gives

$$H_2(M') = H_2(M; \mathbb{Z}[N]) \cong H^1(M; \mathbb{Z}[N]).$$

Our fibration gives rise to a short exact sequence

$$0 \to H^1(N; \mathbb{Z}[N]) \to H^1(M; \mathbb{Z}[N]) \to H^0(N; H^1(M'; \mathbb{Z}[N])) \to \ldots$$

The right hand term is the fixed point set of the (diagonal) action of $N$ on the module $H^1(M'; \mathbb{Z}[N]) \cong H^1(M') \otimes_{\mathbb{Z}} \mathbb{Z}[N]$. By [1, p.69] this is a free $N$-module and since $N$ is an infinite group there are no nontrivial fixed points. Therefore, $H_2(M') \cong H^1(N; \mathbb{Z}[N])$.

**Lemma 8.** *Let $T$ be a finite normal subgroup of a group $N$. Then for all $i \geq 0$ there is an isomorphism $H^i(N; \mathbb{Z}[N]) \cong H^i(N/T; \mathbb{Z}[N/T])$.*

*Proof.* Consider the group extension $1 \to T \to N \to N/T \to 1$ which gives rise to a Serre spectral sequence

$$E_2^{p,q} = H^p(N/T; H^q(T; \mathbb{Z}[N])) \Longrightarrow H^{p+q}(N; \mathbb{Z}[N]).$$

The $E_2^{p,q}$-terms vanish for $q > 0$ because $T$ is a finite group and $\mathbb{Z}[N]$ is a free $T$-module. To prove the lemma it is enough to show that the fixed point set of $T$ on the module $\mathbb{Z}[N]$ is isomorphic (as $N/T$-module) to the group ring $\mathbb{Z}[N/T]$. This isomorphism is the multiplication map

$$\cdot(\sum_{t \in T} t) : \mathbb{Z}[N/T] \to \mathbb{Z}[N]$$

which is a $N/T$-map since the element we are multiplying with lies in the center of $\mathbb{Z}[N]$. An easy calculation shows that this map has image equal to the fixed point set $\mathbb{Z}[N]^T$. It is injective because the composition with the canonical projection $\mathbb{Z}[N] \to \mathbb{Z}[N/T]$ gives multiplication by the order of $T$ which is clearly injective on $\mathbb{Z}[N/T]$. □

Applying this lemma for $i = 1$ to the torsion subgroup $T := \mathrm{Tor}(N)$ of our nilpotent group $N$, we can finish our calculation of $H_2(M')$ because $N/T$ is a Poincaré duality group by Lemma 4. Moreover, the homology groups $H_i(G; \mathbb{Z}[G])$ vanish for any group except for $i = 0$ where one always gets $\mathbb{Z}$. Thus we obtain

$$H_2(M') \cong H_{\mathrm{rank}(N)-1}(N/T; \mathbb{Z}[N/T]) \cong \begin{cases} \mathbb{Z} & \text{if } \mathrm{rank}(N) = 1, \\ 0 & \text{else.} \end{cases}$$

which finishes the proof of Proposition 2.

**Case IIa**: $\mathrm{rank}(N) = 1$, i.e. $N/T$ is infinite cyclic.

We have to prove the following

**Proposition 3.** *In this situation one must have $T = 1$.*

*Proof.* We know that both groups $H_2(M')$ and $H_2(M) \cong H^1(M) \cong H^1(N)$ are infinite cyclic.

**Lemma 9.** *The covering map $M' \to M$ induces multiplication by $|T|$ on $H_2$.*

*Proof.* Consider the commutative diagram

$$
\begin{array}{ccccccc}
H_2(M; \mathbb{Z}[N]) & \xleftarrow[\cong]{\cap[M]} & H^1(M; \mathbb{Z}[N]) & \xleftarrow{\cong} & H^1(N/T; \mathbb{Z}[N/T]) & \xrightarrow[\cong]{\cap[S^1]} & H_0(\mathbb{Z}; \mathbb{Z}[\mathbb{Z}]) \\
\downarrow{\epsilon_N} & & \downarrow{\epsilon_N} & & \downarrow{\epsilon_{N/T}} & & \cong \downarrow{\epsilon_{\mathbb{Z}}} \\
H_2(M; \mathbb{Z}) & \xleftarrow[\cong]{\cap[M]} & H^1(M; \mathbb{Z}) & \xleftarrow{\cdot|T|} & H^1(N/T; \mathbb{Z}) & \xrightarrow[\cong]{\cap[S^1]} & H_0(\mathbb{Z}; \mathbb{Z})
\end{array}
$$

The multiplication by $|T|$ at the bottom line occurs because $\epsilon_N(\sum_{t \in T} t) = |T|$ which has to be taken into account when applying the isomorphism of Lemma 8. The above diagram finishes the proof since the covering map $M' \to M$ induces the same map on $H_2$ as the augmentation map $\epsilon_N : \mathbb{Z}[N] \to \mathbb{Z}$ on the left hand side. □

Looking back to the spectral sequence for the fibration $M' \to M \xrightarrow{f} K(N,1)$ we may now conclude that $H_2(N)$ is cyclic of order $|T|$, generated by the image of $f_*$. Moreover, the commutative diagram

$$
\begin{array}{ccc}
H^1(M) & \xrightarrow[\cong]{\cap [M]} & H_2(M) \\
\cong \Big\uparrow{\scriptstyle f^*} & & \Big\downarrow{\scriptstyle f_*} \\
H^1(N) & \xrightarrow{\cap f_*[M]} & H_2(N)
\end{array}
$$

shows that, with $z \in H^1(N) \cong \mathbb{Z}$ a generator, the map

$$. \cap z : H_3(N) \to H_2(N) \cong \mathbb{Z}/|T|$$

is an epimorphism. Now consider the extension

$$1 \to T \xrightarrow{i} N \to N/T \cong \mathbb{Z} \to 1.$$

The Wang sequence gives a short exact sequence

$$0 \to H_3(T)_{N/T} \to H_3(N) \xrightarrow{c} H_2(T)^{N/T} \to 0.$$

Moreover, the composition $H_3(N) \xrightarrow{c} H_2(T) \xrightarrow{i_*} H_2(N)$ is given by the cap-product with the class $\pm z \in H^1(N)$. But we know that this composition is an epimorphism. In particular, the inclusion map $i$ induces also an epimorphism $i_* : H_2(T) \to H_2(N)$ and thus $H_2(T) \cong H^3(T)$ contains an element of order $|T|$. Then the finite group $T$ has 3-periodic cohomology and is therefore trivial ([1, p.159, exercise 1]). $\qquad\square$

*Proof of Theorem 2.* Again we only have to prove that (4) implies (1). Let $M$ be a closed orientable 3-manifold whose fundamental group allows a maximal *torsionfree* nilpotent quotient $N$. Consider again the fibration

$$M' \to M \to K(N,1)$$

induced from the quotient map $\pi_1 M \to N$. Observe that by assumption $H_1(M')_N$ is a finite group. If one now considers rational homology throughout (e.g. Dwyer's result holds true in this case) then the arguments from the proof of Theorem 1 show that $N$ is trivial, infinite cyclic or a 3-dimensional rational Poincaré duality group. In the last case one sees easily that $N$ is a Heisenberg group. Note that all the details about certain degrees etc. are not necessary here since $N$ is torsionfree to start with. $\qquad\square$

University of California, San Diego, La Jolla, CA, 92093-0112
*E-mail address*: teichner@euclid.ucsd.edu